\documentclass[11pt,letterpaper]{amsart}
\usepackage[T1]{fontenc}
\usepackage{lmodern,amsmath,amssymb,amsthm,mathtools,booktabs,array}
\usepackage{microtype}
\usepackage{xcolor}
\usepackage{etoolbox,needspace}
\usepackage[
    colorlinks=true,
    linkcolor=blue!50!black,
    citecolor=blue!50!black,
    urlcolor=blue!50!black]{hyperref}
\numberwithin{equation}{section}
\allowdisplaybreaks
\newtheorem{theorem}{Theorem}[section]
\newtheorem{proposition}[theorem]{Proposition}
\newtheorem{lemma}[theorem]{Lemma}
\newtheorem{corollary}[theorem]{Corollary}
\theoremstyle{definition}
\newtheorem{definition}[theorem]{Definition}
\theoremstyle{remark}
\newtheorem{remark}[theorem]{Remark}
\BeforeBeginEnvironment{definition}{\Needspace{7\baselineskip}}
\BeforeBeginEnvironment{theorem}{\Needspace{7\baselineskip}}
\BeforeBeginEnvironment{proposition}{\Needspace{6\baselineskip}}
\BeforeBeginEnvironment{corollary}{\Needspace{6\baselineskip}}
\BeforeBeginEnvironment{lemma}{\Needspace{6\baselineskip}}
\BeforeBeginEnvironment{remark}{\Needspace{5\baselineskip}}
\BeforeBeginEnvironment{proof}{\Needspace{4\baselineskip}}
\makeatletter
\newlength{\majorsectionsep}
\patchcmd{\section}{.7\linespacing\@plus\linespacing}{\majorsectionsep}{}{}
\makeatother
\pretocmd{\subsection}{\Needspace{5\baselineskip}}{}{}
\newcommand{\R}{\mathbb R}
\newcommand{\Sph}{\mathbb S}
\newcommand{\Hess}{\operatorname{Hess}}
\newcommand{\tr}{\operatorname{tr}}
\newcommand{\sinc}{\operatorname{sinc}}
\newcommand{\diag}{\operatorname{diag}}
\title[The Lai--Wei construction and minimizing hypercones]
{\small A Generalized Lai--Wei Construction
\\
and New Families of Area-Minimizing Hypercones}
\author{Yongsheng Zhang}
\date{}
\hypersetup{
  pdftitle={A Generalized Lai--Wei Construction and New Families of Area-Minimizing Hypercones},
  pdfauthor={Yongsheng Zhang},
  pdfsubject={Minimal hypersurfaces and area-minimizing hypercones},
  pdfkeywords={area-minimizing hypercone, isoparametric foliation,
        Hsiang--Lawson reduction, minimal hypersurface,
        Lawlor criterion}}

\begin{document}

    \begin{abstract}
        Our results are twofold.
        Inspired by recent work on weighted spheres,
            we extend the Lai--Wei construction from one to several independent isoparametric inputs,
            replacing its profile curve by a hypersurface profile to generate minimal hypersurfaces in spheres.
        Using a curvature criterion,
            we obtain new infinite families of area-minimizing hypercones from the original construction and its extensions.
    \end{abstract}
    \maketitle

\section{Introduction and main results}
\setlength{\majorsectionsep}{1.7\baselineskip}

        Constructing closed minimal hypersurfaces in spheres is a longstanding and attractive problem in differential geometry.
        Classical methods include isoparametric foliations \cite{Cartan},
            min--max theory \cite{Almgren,Pitts},
            reflection \cite{Lawson70},
            symmetry reduction \cite{HL,Hsiang83},
            and gluing \cite{KY}.

        Isoparametric hypersurfaces in spheres have constant principal curvatures;
            their parallel leaves form the regular part of a singular foliation \cite{Munzner}.
        Recently,
            Lai and Wei \cite{LW} extended the product constructions of Carlotto and Schulz \cite{CarlottoSchulz} and Firester and Tsiamis \cite{FiresterTsiamis} to arbitrary isoparametric inputs.
        Their construction combines closed profile curves with the parallel leaves of isoparametric hypersurfaces $P\subset\Sph^n$,
            generating closed embedded minimal hypersurfaces in $\Sph^{n+1}$ with topology $S^1\times P$.

        Inspired by our framework of weighted spheres \cite{MMP},
            we extend the Lai--Wei construction to several independent isoparametric inputs.
        We replace the profile curve by a hypersurface in a coefficient sphere,
            whose coordinates describe the relative sizes and parallel positions of the inputs.
        The leaf volumes define a weighted area on the regular region,
            and a closed embedded stationary profile there generates a closed embedded minimal hypersurface.
        The weight is more general than the monomial weights considered in that framework.
        The main task is therefore to find closed embedded stationary profiles for this weighted area.
        For large multiplicities,
            rescaling near the maximum of the weight yields a generally anisotropic Gaussian density.
        We apply the Gaussian transfer of \cite{MMP} to spherical and rotational profiles in this model,
            obtaining a constellation of closed embedded minimal hypersurfaces in spheres based on isoparametric families of spheres.

    \Needspace{5\baselineskip}
        Closed embedded minimal hypersurfaces in spheres span regular minimal hypercones,
            smooth away from the vertex.
        Regular area-minimizing hypercones are minimal hypercones over embedded closed hypersurfaces 
        with strong global minimizing property, 
        which connect differential geometry with geometric measure theory
        and arise naturally in the regularity theory of the Plateau problem \cite{HS} and in the Bernstein problem \cite{Simon}.
        The classical Simons cones \cite{BDGG,Simons},
            Lawson cones \cite{Lawson72} and Ferus--Karcher cones \cite{FK85} have isoparametric links.
        Recently,
            Cui \cite{Cui} combined Lawlor's curvature criterion \cite{Lawlor} with Howard's rolling theorem \cite{Howard}
             to prove the area minimality for selected Hsiang-type cones,
            answering the existence question by Ilmanen and White \cite{IW} for nonflat,
            topologically trivial hypercones.

       In this paper we establish curvature bounds for the original Lai--Wei hypersurfaces and the families constructed here,
           and obtain new infinite families of regular area-minimizing hypercones.
        Infinitely many of their links are not diffeomorphic to any isoparametric hypersurface of a sphere;
            see Remark~\ref{top:type}.
        We shall focus on isoparametric foliations with $g=2,4$ principal curvatures,
            as they can occur in arbitrarily large dimensions \cite{FKM,Munzner}.

    \begin{theorem}[Several isoparametric inputs]\label{thm:families}
        Let $P_i^{d_i}\subset\Sph^{d_i+1}$,
            $1\leq i\leq k$,
            be regular leaves of isoparametric foliations.
        Assume that $P_2,\ldots,P_k$ have two distinct principal curvatures,
            while $P_1$ has $g_1\in\{2,4\}$ distinct principal curvatures.
        Let $(a_i,b_i)$ be the alternating positive multiplicities,
            and set $D=\sum_i d_i$ and $\mu=\min_i\{a_i,b_i\}$.
        There are an absolute integer $k_0$ and thresholds $\mu_0(k)>0$.
        For each fixed $k\geq2$ and $\mu\geq\mu_0(k)$,
            the following manifolds can be realized as smooth closed embedded minimal hypersurfaces $L^{D+2k-2}\subset\Sph^{D+2k-1}$ spanning area-minimizing hypercones:
        \begin{enumerate}
            \item[\textup{(a)}] $\Sph^{2k-2}\times\prod_{i=1}^kP_i$,
                if $g_1=2$;
            \item[\textup{(b)}] $S^1\times S^{2k-3}\times\prod_{i=1}^kP_i$,
                if $g_1=4$ and $k\geq k_0$.
        \end{enumerate}
    \end{theorem}

        By varying the distribution of dimensions among the inputs in \textup{(a)} and \textup{(b)} while keeping $D$ fixed,
            we obtain links of different homotopy types.

        We now return to the original Lai--Wei construction and obtain the following result.
        For $g\in\{2,4\}$,
            let $a,b$ be the alternating positive multiplicities,
            and set $M=a+b$ and $d=\frac {gM} {2}$.

    \begin{theorem}[Single isoparametric input]\label{thm:single}
        There exist constants $C>0$ and $M_0>0$,
            independent of the input foliation and its multiplicity ratio.
        If $P^d\subset\Sph^{d+1}$ is a regular leaf of an isoparametric foliation with $g\in\{2,4\}$ and $M\geq M_0$,
            then every closed embedded minimal hypersurface $L^{d+1}\subset\Sph^{d+2}$ constructed from $P$ in \cite{LW} spans an area-minimizing hypercone 
            with
      $$
        L\cong S^1\times P,\qquad
            \sup_L|A_L|^2\leq CM.
    $$
    \end{theorem}

        Unlike Theorem~\ref{thm:families},
            this result requires only the sum $a+b$ to be sufficiently large,
            and
            either multiplicity may remain fixed.

    \Needspace{7\baselineskip}
        Finally,
            the same viewpoint applies to a whole sphere factor $\Sph^d$ and a freely varying point in $\R^3$.
        The sphere contributes one radial coordinate $u>0$,
            and the point contributes three free coordinates.
        Thus we consider  a coefficient hemisphere of $\Sph^3$ with weight $u^d$.
        We adapt the doubling construction of Kapouleas and McGrath~\cite{KM} to the corresponding conformal metrics,
        and find closed embedded stationary surface profiles of genus $h$ that generate minimal hypersurfaces in $\Sph^{d+3}$.
        Let $\Sigma_h$ denote a closed connected orientable surface of genus $h$.

    \begin{theorem}\label{thm:genus}
        For every sufficiently large odd integer $h$,
            there are constants $C_h>0$ and $d_0(h)$ such that,
            for every integer $d\geq d_0(h)$,
            there exists a smooth closed embedded minimal hypersurface $L_{d,h}^{d+2}\subset\Sph^{d+3}$ spanning an area-minimizing hypercone with
    $$
        L_{d,h}\cong\Sigma_h\times\Sph^d,\qquad
            \sup_{L_{d,h}}|A_{L_{d,h}}|^2\leq C_h d.
    $$
    \end{theorem}

        For $d\geq2$,
            the links in Theorem~\ref{thm:genus} have first Betti number $b_1(L_{d,h})=2h$.

        The paper is organized as follows.
        In \S\ref{sec:construction} we present the generalized Lai--Wei construction.
        In \S\ref{sec:gaussian}--\S\ref{sec:genus} 
        the closed profiles for Theorems~\ref{thm:families} and~\ref{thm:genus} 
           will be constructed by employing Gaussian transfer and doubling.
        Finally, we prove the area-minimizing assertions in \S\ref{sec:criterion}.

    \Needspace{10\baselineskip}
\section{The generalized Lai--Wei construction}\label{sec:construction}

        We derive the weight on the coefficient sphere from the tangential scaling of the input leaves.
        With this weight,
            minimality of the resulting hypersurface is equivalent to stationarity of its profile.
        The input leaves need not be minimal.

        In \cite{LW},
            the authors start with a regular isoparametric hypersurface $P^d\subset\Sph^{d+1}$ with unit normal $\nu$ and consider the immersion $(p,s)\mapsto(u(s)p+v(s)\nu(p),w(s))$,
            generated by a curve $(u,v,w):I\to\Sph^2$.
        For a principal vector $X$ with curvature $\kappa_j$,
            our convention $A=-\nabla\nu$ implies $D_X\nu=-\kappa_jX$,
            and hence $D_X(up+v\nu)=(u-\kappa_jv)X$.
        A principal space of dimension $m_j$ therefore contributes $|u-\kappa_jv|^{m_j}$ to the fiber volume density,
            so $W_1=\prod_j|u-\kappa_jv|^{m_j}$.
        Thus a profile curve generates a minimal hypersurface precisely when it is a geodesic of $W_1^2g_{\Sph^2}$.
        Their shooting and reflection argument yields a closed embedded profile.

        For $k\geq2$ independent inputs,
            the position and normal coefficients lie in $\Sph^{2k-1}$.
        We call a submanifold $\Gamma$ of this sphere a coefficient profile and refer to the resulting submanifold as its lift.
        Hypersurface links require $\dim\Gamma=2k-2$;
            with two inputs the profile is a surface in $\Sph^3$.
        Fix regular isoparametric hypersurfaces $P_i^{d_i}\subset\Sph^{d_i+1}$ with unit normals $\nu_i$,
            and set $D=\sum_i d_i$ and $p=2k-2$.
        Using the decomposition $\R^{D+2k}=\bigoplus_i\R^{d_i+2}$,
            we define
    \begin{equation}\label{mk:lift}
        \mathcal E(p_1,\ldots,p_k,z)=(u_ip_i+v_i\nu_i(p_i))_{i=1}^k,
            \quad z=(u_i,v_i)_{i=1}^k\in\Sph^{2k-1}.
    \end{equation}
    \Needspace{8\baselineskip}
        Let $I_i$ be the maximal regular interval of the $i$th parallel family containing zero.
        The coefficient chamber $Q$ consists of $z=(\rho_i\cos t_i,\rho_i\sin t_i)_i$ with $\rho_i>0$,
            $\sum_i\rho_i^2=1$ and $t_i\in I_i$.
        For the principal curvatures $\kappa_{ij}$ and multiplicities $m_{ij}$ choose signs so that
    $$
        g_Q=g_{\Sph^{2k-1}}|_Q,\qquad
            \ell_{ij}=\pm(u_i-\kappa_{ij}v_i)>0,
            \qquad W=\prod_{i,j}\ell_{ij}^{m_{ij}}.
    $$
        The factors calculated for one input multiply across the mutually orthogonal blocks.
        The isoparametric assumption makes each $\kappa_{ij}$ constant on $P_i$,
            so the fiber volume density depends only on the coefficient point $z$.
        Here and in the transfer results $k$ is fixed.
        The separate large-$k$ limit needed for Theorem~\ref{thm:families}\textup{(b)} is treated in \S\ref{sec:ring}.

    \begin{lemma}[Weighted reduction and curvature]\label{mk:reduction}
        The map \eqref{mk:lift} is a diffeomorphism from $\prod_iP_i\times Q$ onto its image.
        Every smooth closed embedded hypersurface $\Gamma^p\subset Q$ lifts to a closed embedded hypersurface $L^{D+p}\subset\Sph^{D+p+1}$ with $L\cong\Gamma\times\prod_iP_i$.
        Under this identification,
            the area element is $dA_L=W\,dA_{g_Q}\prod_i dV_{P_i}$.
        With trace mean curvature $H=\tr A$ and a unit normal $\eta$ to $\Gamma$ in $(Q,g_Q)$,
    \begin{equation}\label{mk:hessian-curvature}
        H_L=H_\Gamma-\eta\log W,\quad
            |A_L|^2=|A_\Gamma|_{g_Q}^2-D-
            \Hess_{g_Q}\log W(\eta,\eta).
    \end{equation}
        Thus $L$ is minimal precisely when $\Gamma$ is stationary for $\int_\Gamma W\,dA_{g_Q}$,
            or equivalently minimal in $W^{\frac {2} {p}}g_Q$.
    \end{lemma}

    \begin{proof}
        The first assertion follows from regular normal coordinates.
        The area formula follows from the above scalings and orthogonality of coefficient and leaf directions.
        We extend the principal-curvature calculation in \cite[\S2]{LW} to the independent input blocks.
        For $\eta=(\eta_{u_i},\eta_{v_i})_i$,
            the lifted normal is $(\eta_{u_i}p_i+\eta_{v_i}\nu_i)_i$.
        The shape operator has horizontal block $A_\Gamma$,
            zero mixed block and fiber eigenvalues $-\eta\log\ell_{ij}$ of multiplicities $m_{ij}$.
        Taking the trace and squared norm,
            and using $\Hess\log\ell_{ij}=-g_Q-d\log\ell_{ij}\otimes d\log\ell_{ij}$ for the restricted linear forms,
            we obtain \eqref{mk:hessian-curvature}.
        The trace identity and the conformal area factor prove the stationarity equivalences;
            cf.~\cite{HL,PT}.
    \end{proof}

    \begin{remark}[Topology]\label{top:type}
        Let $\beta(X)=\sum_j\dim H^j(X;\mathbb Z_2)$ be the total Betti number with coefficients in $\mathbb Z_2$.
        M\"unzner proved $\beta(P)=2g\leq12$ for closed connected embedded isoparametric hypersurfaces of spheres \cite{Munzner}.
        For the families in the main theorems,
            the product topologies and the K\"unneth formula give
    \begingroup
        \renewcommand{\arraystretch}{1.2}
    $$
        \begin{array}{c|cccc}
            \text{Theorem}&\ref{thm:families}\textup{(a)}
                &\ref{thm:families}\textup{(b)}
                &\ref{thm:single}\ (g=4)&\ref{thm:genus}\\ \hline
            \beta(L)&2^{2k+1}&2^{2k+3}&16&4(h+1).
        \end{array}
    $$
    \endgroup
        Each of these links has total Betti number greater than $12$,
            so it cannot be diffeomorphic to an isoparametric hypersurface of a sphere.

        The complement of the two focal submanifolds of each input is diffeomorphic to $P\times(-1,1)$ \cite{Munzner}.
        If $\mu\geq2$,
            both focal submanifolds have codimension at least $3$.
        Every loop in their complement bounds a disk in the sphere.
        Since the disk is two-dimensional,
            relative transversality allows it to avoid both focal submanifolds while fixing its boundary.
        Hence $P$ is simply connected.
        After choosing $\mu_0(k)\geq2$ and $k_0\geq3$ in Theorem~\ref{thm:families},
            the links in \textup{(a)} and \textup{(b)} satisfy $\pi_1(L)=0$ and $\pi_1(L)\cong\mathbb Z$,
            respectively.
    \end{remark}

\section{Gaussian profiles for the extended construction}\label{sec:gaussian}

        Near the maximum of the weight,
            we identify the Gaussian limit and verify the chamber condition needed for the transfer theorem of \cite{MMP}.

\subsection{Uniform anisotropic Gaussian limit}

        We fix $g_i\in\{2,4\}$ and let $a_i,b_i$ be the alternating multiplicities of the $i$th input.
        Set $q_i=\frac {g_i} {2}$,
            $d_i=q_i(a_i+b_i)$ and $\mu=\min_i\{a_i,b_i\}$.
        Translate the angular coordinates so that $0<t_i<\frac {\pi} {g_i}$.
        The sine product formula gives
    \begin{equation}
        W=c\prod_i\rho_i^{d_i}
            \sin^{a_i}(q_it_i)\cos^{b_i}(q_it_i),\qquad c>0,
    \end{equation}
        and $g_Q=\sum_i d\rho_i^2+\sum_i\rho_i^2dt_i^2$,
            subject to $\sum_i\rho_i^2=1$.
        The unique maximum $z_*$ of $W$ satisfies $\rho_{i*}^2=\frac {d_i} {D}$ and $\sin^2(q_it_{i*})=\frac {a_i} {a_i+b_i}$.
        At $z_*$,
            let $e_i=(\cos t_{i*},\sin t_{i*})$ and $f_i=(-\sin t_{i*},\cos t_{i*})$ be supported in the $i$th coefficient plane.
        With $\rho_*=(\rho_{1*},\ldots,\rho_{k*})$,
            choose $U\in\R^{k\times(k-1)}$ such that $U^TU=I$ and $U^T\rho_*=0$.
        We obtain an orthonormal frame by $E(x,y)=\sum_i((Ux)_ie_i+y_if_i)$,
            where $(x,y)\in\R^{k-1}\oplus\R^k$.
        In this frame,
            we have
    \begin{equation}\label{mk:hessian}
        E^*(\Hess_{g_Q}\log W(z_*))=-D\mathsf A,
            \qquad
            \mathsf A=\diag(2I_{k-1},g_1,\ldots,g_k).
    \end{equation}
        The ambient Hessian has entries $-D$ and $(1-g_i)D$ in the directions $e_i,f_i$.
        The radial derivative of $\log W$ is $D$,
            so restriction to the unit sphere subtracts $Dg_Q$.
        Thus the Gaussian anisotropy depends only on $(g_i)$,
            while the spatial rescaling below is isotropic.

    \begin{lemma}[Uniform Gaussian limit]\label{mk:gaussian-limit}
        For fixed $R$ and each fixed derivative order $s\geq0$,
            define
    $$
        \Phi_D(X)=\exp_{z_*}(\frac {EX} {\sqrt D}),\qquad
            g_D=D\Phi_D^*g_Q,\qquad
            \psi_D=\log\frac {W\circ\Phi_D} {W(z_*)}.
    $$
        Then,
            uniformly over all multiplicity ratios and all choices of input foliations,
    \begin{equation}\label{mk:limit}
        g_D=g_{\mathrm E}+O_{C^s(B_R)}(D^{-1}),\qquad
            \psi_D=-\frac {1} {2}\langle\mathsf AX,X\rangle
            +O_{C^s(B_R)}(\mu^{-\frac {1} {2}}).
    \end{equation}
        For sufficiently large $\mu$,
            we have $\Phi_D(B_R)\subset Q$.
    \end{lemma}

    \begin{proof}
        We combine the indices into $j$ and set $b_j=\frac {E^*\nabla\ell_j(z_*)} {\sqrt D\,\ell_j(z_*)}$ and $T=\tr(\mathsf A-I)=\sum_i g_i-1$.
        By the Hessian identity for linear forms in Lemma~\ref{mk:reduction} and \eqref{mk:hessian},
            we have
    $$
        \sum_jm_jb_j=0,\quad
            \sum_jm_jb_j\otimes b_j=\mathsf A-I,\quad
            \max_j|b_j|\leq\sqrt{\frac {T} {\mu}}.
    $$
        With $\sinc t=\frac {\sin t} {t}$,
            $\sinc0=1$ and $\zeta=\frac {|X|} {\sqrt D}$,
            we have the exact identity
    $$
        \frac {\ell_j(\Phi_D(X))} {\ell_j(z_*)}
            =\cos\zeta+\sinc(\zeta)b_j\cdot X.
    $$
        For fixed $R$,
            these ratios are positive for sufficiently large $\mu$,
            so $\Phi_D(B_R)\subset Q$.
        Since $\sum_jm_j=D$ and $\sum_jm_j|b_j|^r\leq T(\frac {T} {\mu})^{\frac {r-2} {2}}$ for $r\geq2$,
            Taylor expansion and the moment identities prove the density estimate in \eqref{mk:limit}.
        The metric estimate follows from the scaled round exponential map.
    \end{proof}

\subsection{Gaussian transfer and curvature estimates}

        We apply the anisotropic Gaussian transfer of \cite[Lemma~4.2 and the remark on p.~15]{MMP},
            and then compute the lifted curvature.
        We use the trivial symmetry group and require the full Jacobi kernel to vanish,
            since the finite weights need not preserve the symmetries of the limit.

    \begin{theorem}[Gaussian transfer]\label{thm:gaussian-seed}
        With the inputs and notation above,
            fix a smooth closed embedded two-sided hypersurface $\Gamma_\infty^p\subset\R^{p+1}$,
            stationary for $\int_{\Gamma_\infty}e^{-f}\,dA$,
            where $f(X)=\frac {\langle\mathsf AX,X\rangle} {2}$.
        Assume its Gaussian Jacobi kernel vanishes among all smooth normal variations.
        For sufficiently large $\mu$,
            there is a smooth closed embedded minimal link $L^{D+p}\subset\Sph^{D+p+1}$ with $L\cong\Gamma_\infty\times\prod_iP_i$.
        Under the rescaling of Lemma~\ref{mk:gaussian-limit} and pullback to $\Gamma_\infty$,
            we have uniform convergence $\frac {|A_L|^2} {D}\to|A_{\Gamma_\infty}|^2+\Hess f(\nu,\nu)-1$,
            where $\nu$ is a unit normal to $\Gamma_\infty$.
        For the fixed Gaussian seed,
            $k$ and $(g_i)$,
            the threshold and convergence are uniform over all multiplicity ratios and choices of input foliations.
    \end{theorem}

    \begin{proof}
        Fix $0<\alpha<1$ and a tubular neighborhood of $\Gamma_\infty$.
        The $C^{3,\alpha}$ convergence required by \cite[Lemma~4.2]{MMP} follows from Lemma~\ref{mk:gaussian-limit} with $s=4$.
        We obtain a stationary normal graph $\widehat\Gamma_D=\{q+u_D(q)\nu(q):q\in\Gamma_\infty\}$,
            with $\|u_D\|_{C^2}=O(\mu^{-\frac {1} {2}})$.
        Lemma~\ref{mk:reduction} lifts $\Gamma_D=\Phi_D(\widehat\Gamma_D)$ to the asserted minimal link.

        Pulling back \eqref{mk:hessian-curvature} by the normal graph and $\Phi_D$,
            we obtain $\frac {|A_L|^2} {D}=|A_{\widehat\Gamma_D}|_{g_D}^2-1-\Hess_{g_D}\psi_D(\nu_D,\nu_D)$,
            where $\nu_D$ is a $g_D$-unit normal.
        The asserted limit follows from Lemma~\ref{mk:gaussian-limit},
            with uniform error $O(\mu^{-\frac {1} {2}})$.
    \end{proof}

        When all $g_i=2$,
            we have $\mathsf A=2I$.
        The round Gaussian hypersurface $\Gamma_\infty=\Sph^{2k-2}(\sqrt{k-1})\subset\R^{2k-1}$ has trivial Jacobi kernel by \cite[Lemma~4.4]{MMP}.
        We obtain the family in Theorem~\ref{thm:families}\textup{(a)} by Theorem~\ref{thm:gaussian-seed},
            with $\frac {|A_L|^2} {D}\to3$ uniformly since $|A_{\Gamma_\infty}|^2=2$.

\subsection{Rotational Gaussian hypersurfaces}\label{sec:ring}

        To obtain the family with mixed inputs in Theorem~\ref{thm:families}\textup{(b)},
            we need a closed fully nondegenerate hypersurface for an anisotropic Gaussian density.
        We construct it by rotating a closed curve that approaches a circle after a translation of the radial coordinate,
            as the dimension of the spherical factor tends to infinity.
        For one $g=4$ input and $k-1$ inputs with $g=2$,
            set $N=2k-3$.
        The Gaussian density is $e^{-2x^2-|Y|^2}$ on $\R\times\R^{N+1}$,
            with $\mathsf A=\diag(4,2I_{N+1})$ by Lemma~\ref{mk:gaussian-limit}.
        After dilating all coordinates by $\sqrt2$,
            the equivalent density is $e^{-x^2-\frac {|Y|^2} {2}}$ on $\R\times\R^{N+1}$.
        An $O(N+1)$-invariant hypersurface is generated by a curve $\gamma=(x,y)$ in the half-plane $y=|Y|>0$.
        Apart from the area of $S^N$,
            its weighted area is
    \begin{equation}\label{rg:functional}
        \mathcal F_N(\gamma)=\int_\gamma y^N e^{-x^2-\frac {y^2} {2}}\,ds.
    \end{equation}

        Since $y^Ne^{-\frac {y^2} {2}}$ is maximal at $y=\sqrt N$,
            we take $\varepsilon=N^{-\frac {1} {2}}$ and $y=\varepsilon^{-1}+z$.
        After discarding a constant factor,
            \eqref{rg:functional} becomes $\int e^{-F_\varepsilon}\,ds$ with $F_\varepsilon(x,z)=x^2+\frac {z^2} {2} +\frac {\varepsilon z-\log(1+\varepsilon z)} {\varepsilon^2}$.
        Centering at the maximum cancels the linear term in $z$,
            leaving an isotropic quadratic term:
    \begin{equation}
        F_\varepsilon(x,z)=x^2+z^2-\frac {\varepsilon z^3} {3}
            +\frac {\varepsilon^2z^4} {4}+O(\varepsilon^3).
    \end{equation}
        The potential $F_\varepsilon$ extends smoothly to $\varepsilon=0$ on every fixed compact set.

        Stationary meridians satisfy $\kappa=\nabla F_\varepsilon\cdot\nu$,
            where $\nu=(\nu_x,\nu_z)$ is the outward unit normal and $\kappa=-H_\gamma$ is positive on a circle.
        At $\varepsilon=0$,
            the circle of radius $c=\frac {1} {\sqrt2}$ is fully nondegenerate by \cite[Lemma~4.4]{MMP}.
        Applying \cite[Lemma~4.2]{MMP},
            we obtain a unique smooth family of stationary radial graphs $(x,z)=r_\varepsilon(\theta)(\cos\theta,\sin\theta)$ near this circle.
        For small positive $\varepsilon$,
            the meridian is embedded and $y=\varepsilon^{-1}+z>0$.
        Rotating this meridian,
            we obtain a closed embedded hypersurface $\mathcal T_N\cong S^1\times S^N\subset\R^{N+2}$.
        The meridian equation is precisely its weighted mean curvature equation,
            so $\mathcal T_N$ is stationary under all normal variations.

        To apply Theorem~\ref{thm:gaussian-seed},
            we decompose normal functions into spherical harmonics on the $S^N$ fibers.
        With meridian arclength $\sigma$,
            the degree-$\ell$ Jacobi operator is
    $$
        \begin{aligned}
            J_{0,\varepsilon}
                &=\partial_\sigma^2-(\partial_\sigma F_\varepsilon)\partial_\sigma
                    +\kappa^2+\Hess F_\varepsilon(\nu,\nu),\\
            J_{\ell,\varepsilon}
                &=J_{0,\varepsilon}
                    -\frac {\ell[1+\varepsilon^2(\ell-1)]} {(1+\varepsilon z)^2}.
        \end{aligned}
    $$
        Indeed,
            $|A_{\mathcal T_N}|^2=\kappa^2+\frac {N\nu_z^2} {y^2}$ and $\Hess F_\varepsilon=\diag(2,1+\frac {N} {y^2})$.
        These operators are self-adjoint in $L^2(e^{-F_\varepsilon}\,d\sigma)$ and have smoothly varying coefficients on the fixed circle parametrized by $\theta$.

        At the circle,
            $J_{\ell,0}=2\partial_\theta^2+4-\ell$.
        Its eigenvalue at meridian Fourier frequency $j$ is $4-\ell-2j^2$.
        Zero occurs only for $(\ell,j)=(4,0)$ and $(2,1)$,
            with meridian eigenfunctions $1$ and $\cos\theta,\sin\theta$,
            respectively.

    \begin{proposition}[Jacobi nondegeneracy]\label{rg:nondegenerate}
        For every sufficiently large integer $N$,
            the rotational Gaussian hypersurface $\mathcal T_N$ is fully Jacobi nondegenerate.
    \end{proposition}

    \begin{proof}
        We use translations and dilation to determine the signs of the three eigenvalues tending to zero.
        By uniqueness,
            reflection in $x$ preserves the meridian and separates the two $\ell=2$ modes.
        The limiting zero eigenvalue is simple in each relevant reflection class.
        Set $s=\sin\theta$,
            $t=\cos\theta$ and $q=x\nu_x+z\nu_z$.
        We linearize $\kappa=\nabla F_\varepsilon\cdot\nu$ at the circle to obtain $r_\varepsilon=c+\varepsilon u+O(\varepsilon^2)$,
            where $\partial_\theta^2u+2u=\frac {s^3} {4}$ and $u=\frac {6s-s^3} {28}$.
        Differentiating the meridian equation under translations and dilation,
            we have
    $$
        J_{0,\varepsilon}\nu_x=2\nu_x,\qquad
            J_{0,\varepsilon}\nu_z=\bigl[1+(1+\varepsilon z)^{-2}\bigr]\nu_z,
    $$
    $$
        J_{0,\varepsilon}q=4q-
            \frac {\varepsilon z^2(3+2\varepsilon z)} {(1+\varepsilon z)^2}\nu_z.
    $$
    \Needspace{10\baselineskip}
        Let $\langle\cdot,\cdot\rangle_\varepsilon$ denote the inner product in $L^2(e^{-F_\varepsilon}\,d\sigma)$.
        Using these identities and the product rule,
            we obtain $\psi_g=g+O(\varepsilon)$ and $J_{\ell,\varepsilon}\psi_g=O(\varepsilon^2)$ uniformly on the circle for the functions below,
            with
    $$
        \langle g,J_{\ell,\varepsilon}\psi_g\rangle_\varepsilon
            =\varepsilon^2c_g\langle g,g\rangle_\varepsilon
            +O(\varepsilon^3).
    $$
    \begin{center}\small
        \setlength{\tabcolsep}{4pt}
        \renewcommand{\arraystretch}{1.5}
        \begin{tabular}{@{}c c c@{}}
            \toprule
            $(\ell,g)$ & $\psi_g$ & $c_g$\\
            \midrule
            $(4,1)$
                & $\frac {q} {c}+\varepsilon(3z-\frac {z^3} {3})$
                & $-\frac {33} {4}$\\
            $(2,t)$
                & $(1+\varepsilon z)^{\frac {2} {3}}\nu_x$
                & $-\frac {13} {7}$\\
            $(2,s)$
                & $(1+\varepsilon z)^{\frac {1} {3}}(\nu_z-\frac {2\varepsilon q} {3})$
                & $-\frac {93} {28}$\\
            \bottomrule
        \end{tabular}
    \end{center}
        Let $e_{g,\varepsilon}=g+O(\varepsilon)$ be the corresponding eigenfunction with eigenvalue $\lambda_g$.
        By self-adjointness,
            $\lambda_g\langle e_{g,\varepsilon},\psi_g\rangle_\varepsilon
            =\langle e_{g,\varepsilon},J_{\ell,\varepsilon}\psi_g\rangle_\varepsilon$,
            hence $\lambda_g=\varepsilon^2c_g+O(\varepsilon^3)<0$.

        For $0\leq\ell\leq3$,
            uniform ellipticity and spectral continuity exclude zero from the remaining spectrum.
        The eigenfunction of $J_{4,\varepsilon}$ continuing $1$ stays positive,
            so its negative eigenvalue is the largest.
        The angular potential increases with $\ell$,
            hence $J_{\ell,\varepsilon}\leq J_{4,\varepsilon}<0$ as quadratic forms for $\ell\geq4$.
    \end{proof}

        For every fixed $k\geq k_0$,
            the seed $\Gamma_\infty=\frac {1} {\sqrt2}\mathcal T_{2k-3}$ in the original coefficient coordinates is fully nondegenerate by Proposition~\ref{rg:nondegenerate}.
        Applying Theorem~\ref{thm:gaussian-seed},
            we obtain the minimal links in Theorem~\ref{thm:families}\textup{(b)}.
        The FKM families \cite{FKM} supply $g=4$ inputs with both multiplicities arbitrarily large,
            so this case occurs in infinitely many dimensions.

\section{Surface profiles of unbounded genus}
\label{sec:genus}

        In this section we will apply the doubling theorem of Kapouleas and McGrath~\cite{KM} to a stationary surface in the coefficient hemisphere,
            while keeping the number of bridges fixed as the input dimension increases.

\subsection{The weighted coefficient hemisphere}

        For $Q=\{(u,w)\in\Sph^3:u>0,\ w\in\R^3\}$ and $W_d=u^d$,
            an embedded surface $\Gamma\subset Q$ has the lift
    \begin{equation}\label{hg:lift}
        F:\Gamma\times\Sph^d\longrightarrow\Sph^{d+3},\qquad
            F(q,p)=(u(q)p,w(q)).
    \end{equation}
        As in Lemma~\ref{mk:reduction},
            $F$ is an embedding.
        For $L=F(\Gamma\times\Sph^d)$ and a unit normal $\eta$ to $\Gamma$ in the coefficient sphere,
            we have
    \begin{equation}\label{hg:shape}
        H_L=H_\Gamma-d\,\eta\log u,\qquad
            |A_L|^2=|A_\Gamma|^2+d(\eta\log u)^2.
    \end{equation}
        Thus $L$ is minimal precisely when $\Gamma$ is stationary for $\int_\Gamma W_d\,dA$,
            equivalently minimal for $\bar g_d=W_dg_{\Sph^3}$.
        We start from the stationary latitude $u^2=\frac {d} {d+2}$,
            whose lift is a minimal Clifford hypersurface.

        Set $a_d^2=\frac {2} {d+2}$ and $\epsilon=\frac {2} {d}$.
        In the coordinates $(u,w)=(\sqrt{1-a_d^2|x|^2},a_dx)$,
            the seed is $\Sigma=\{|x|=1\}$ and the pullback of $\bar g_d$ divided by $a_d^2(1-a_d^2)^{\frac {d} {2}}$ is
    \begin{equation}
        \widehat g_\epsilon=P_\epsilon^{\frac {1} {\epsilon}}
            \left(g_{\mathrm E}+\frac {\epsilon} {P_\epsilon}(x\cdot dx)^2\right),
            \qquad P_\epsilon=1+\epsilon(1-|x|^2).
    \end{equation}
        The expansion of $\epsilon^{-1}\log P_\epsilon$ shows that the metrics extend smoothly to $\widehat g_0=e^{1-|x|^2}g_{\mathrm E}$ on a fixed compact neighborhood of $\Sigma$.
        At $r=|x|=1$,
            the orbit radius $b(r)=rP_\epsilon^{\frac {1} {2\epsilon}}$ satisfies $b=1$,
            $\partial_\tau b=0$ and $\partial_\tau^2b=-2$,
            where $\tau$ is radial arclength.
        Thus $\Sigma$ is totally geodesic with induced unit round metric $\sigma$ and $\operatorname{Ric}(\nu,\nu)=4$ for every small $\epsilon\geq0$.
        With $\Delta_\sigma=\operatorname{div}_\sigma\nabla_\sigma$,
            the Jacobi operator $\Delta_\sigma+4$ has nonzero eigenvalues $4-\ell(\ell+1)$.

        The radial metric has the $O(2)\times\mathbb Z_2$ symmetry required by the doubling construction;
            the equatorial reflection reverses the orientation of the two-sided $\Sigma$.
        Its Jacobi equation agrees with that of the spherical shrinker in \cite[Theorem~10.5]{KM},
            so we use their two-circle configuration with $m$ bridge centers on each circle.

\subsection{The doubling construction}

    \begin{lemma}\label{hg:fixed}
        There is an integer $m_0$ such that,
            for each integer $m\geq m_0$,
            there are $\epsilon_m>0$ and $C_m<\infty$ with the following property.
        For every $0\leq\epsilon\leq\epsilon_m$,
            $\widehat g_\epsilon$ has a closed embedded minimal surface $M_{\epsilon,m}$ of genus $2m-1$ doubling $\Sigma$,
            with $m$ bridges near each of two symmetric parallel circles.
        These surfaces lie in a common compact neighborhood of $\Sigma$ and satisfy $\sup_{M_{\epsilon,m}}|A|^2_{\widehat g_\epsilon}\leq C_m$.
    \end{lemma}

    \begin{proof}
        Choose $0<\beta'<\beta<1$ and the auxiliary constants of the Kapouleas--McGrath construction valid at both exponents.
        At $\epsilon=0$,
            we use their two-circle family with aligned bridge centers \cite[Theorems~9.39--9.40 and~10.5]{KM}.
        Let $\tau_*$ be the neck scale at the central parameter and $\alpha>0$ their small auxiliary exponent.
        The neck scales are comparable to $\tau_*$ over the parameter box,
            and $\tau_*\to0$ as $m\to\infty$ \cite[Lemma~9.36]{KM}.
        In the fixed point map \cite[(5.10)--(5.11)]{KM},
            the prescribed dislocation cancels and the leading matching term lies in half its parameter box.
        The graph,
            dislocation and remaining matching errors,
            normalized by their defining radii,
            are respectively $O(\tau_*^{1-\frac {5\alpha} {4}})$,
            $O(\tau_*^{\frac {\alpha} {20}})$ and $O(\tau_*^{\frac {\alpha} {4}})$.
        These bounds follow from items \textup{(2)} and \textup{(4)} in the proof of \cite[Theorem~5.7]{KM},
            where the polynomial losses have already been absorbed.
        The fixed sphere's Jacobi inverse controls the obstruction values and first derivatives uniformly in $m$;
            the parameter normalizations in \cite[Definitions~5.3 and~9.32]{KM} give the stated exponents.
        We therefore choose $m_0$ so that all normalized output bounds are at most $\frac {3} {4}$ and the graph correction satisfies a strict smallness bound ensuring embeddedness.

        Fix $m\geq m_0$ and retain the initial embeddings,
            cutoffs,
            substitute kernel extensions,
            intrinsic matching data and finite symmetry group at $\epsilon=0$.
        The neck scales have a positive minimum on the compact parameter set,
            so the identifications in \cite[Lemma~5.5]{KM} make the weighted norms uniformly equivalent to ordinary H\"older norms on a fixed reference surface.
        Let $B_m$ be the product of the prescribed closed invariant weighted $C^{2,\beta}$ graph ball and the parameter boxes.
        With the $C^{2,\beta'}$ topology on the graph factor,
            this set is compact and convex.

        Using $\widehat g_\epsilon$-normal graphs over these fixed embeddings,
            let $L_\epsilon$ be the linearized mean curvature operator and $J_0$ the retained substitute kernel extension.
        For a graph correction $u$ and substitute kernel data $w^\pm$,
            the augmented operators $\mathcal A_\epsilon(u,w^+,w^-)=L_\epsilon u-J_0(w^+,w^-)$ converge to $\mathcal A_0$ in operator norm,
            uniformly over the parameter set at both H\"older exponents.
        For each initial surface,
            the normalized solves in \cite[Proposition~4.18 and Corollary~4.21]{KM} provide the same right inverse bounded at both exponents.
        We obtain corresponding right inverses for small $\epsilon$ by perturbation using a Neumann series.
        The nonlinear operators converge uniformly in $C^{0,\beta}$ on the whole bounded $C^{2,\beta}$ graph ball by smooth convergence of the mean curvature coefficients.
        Using the same extension $J_0$ in both augmented solves,
            we define $T_\epsilon$ by \cite[(5.11)]{KM}.
        These maps are continuous in the weaker topology and satisfy
    $$
        \sup_{z\in B_m}\|T_\epsilon(z)-T_0(z)\|_m\longrightarrow0,
    $$
        where $\|\cdot\|_m$ is the maximum of the weighted $C^{2,\beta}$ graph norm and the parameter norms,
            each normalized by its defining radius.
        For $0\leq\epsilon\leq\epsilon_m$ with $\epsilon_m$ sufficiently small,
            the strict margins persist,
            and a fixed point exists by Schauder's theorem.
        Its graph equation cancels the linear correction in item \textup{(5)} of the proof of \cite[Theorem~5.7]{KM},
            while its parameter equations set the symmetric and antisymmetric matching errors to zero.
        Together,
            these identities force the mean curvature to vanish,
            so the resulting surface is an embedded minimal doubling of genus $2m-1$.
        The positive neck scales and strict graph bound ensure a common compact neighborhood,
            uniform ordinary $C^2$ bounds and uniformly nondegenerate induced metrics,
            hence the asserted curvature bound.
    \end{proof}

        Fix $h=2m-1$ with $m\geq m_0$.
        For every sufficiently large $d$,
            the coordinate map above sends $M_{\frac {2} {d},m}$ to an embedded surface $\Gamma_{d,h}\subset Q$ stationary for weighted area.
        Its lift is a minimal hypersurface $L_{d,h}\cong\Sigma_h\times\Sph^d$.

\section{Area minimization of the cones}
\label{sec:criterion}

        For a closed oriented link $L^n\subset\Sph^{n+1}$,
            let $C(L)=\{rx:r\geq0,\ x\in L\}\subset\R^{n+2}$.
        We call $C(L)$ area minimizing if its truncation to the unit ball minimizes mass among compactly supported integral $(n+1)$-currents in $\R^{n+2}$ with boundary $L$.
        By dilation and comparison in a larger ball,
            this is equivalent to local mass minimization among integral currents in $\R^{n+2}$.

        By Lawlor's curvature criterion \cite{Lawlor} in the form stated by Cui \cite[Proposition~1]{Cui},
            which incorporates Howard's rolling theorem \cite{Howard},
            a smooth closed embedded minimal hypersurface $L^n\subset\Sph^{n+1}$ of the unit round sphere spans an area-minimizing hypercone provided
    \begin{equation}\label{eq:criterion}
        n+1\geq12,\qquad \sup_L|A_L|^2\leq\frac {(n+1)^2} {8}.
    \end{equation}
        In particular,
            a bound $\sup_L|A_L|^2\leq C(n+1)$ suffices once $n+1\geq\max\{12,8C\}$.

        For the original Lai--Wei family,
            we need a bound linear in dimension even when one multiplicity stays fixed.

\subsection{Uniform curvature bounds for the original links}
\label{sf:setup}

        We prove Theorem~\ref{thm:single} for $q=\frac {g} {2}\in\{1,2\}$,
            $M=a+b$,
            $d=qM$ and $n=d+1$.
        The profile estimates hold for all integers $a,b\geq1$;
            an actual foliation is needed only for the geometric lift.

        Let $r\in(0,\pi)$ be the suspension coordinate,
            and let $\theta\in(0,\frac {\pi} {g})$ parametrize the regular leaves of the given isoparametric foliation.
        Use $\xi=q\log\tan(\frac {r} {2})$,
            $v=q\theta\in(0,\frac {\pi} {2})$ and $v_* =\arctan\sqrt{\frac {a} {b}}$.
        We set $\mathfrak a(\xi)=\frac {n} {q}\tanh\frac {\xi} {q}$ and
    $$
        G(v)=a\cot v-b\tan v,\qquad
            Q(v)=a\csc^2v+b\sec^2v=-G'(v).
    $$
        The reduced metric $g_Q=dr^2+\sin^2r\,d\theta^2$ becomes $q^{-2}\operatorname{sech}^2(\frac {\xi} {q})(d\xi^2+dv^2)$.
        After a constant normalization,
            the volume factor of the leaves is $W=\sin^d r\,\sin^a v\cos^b v$.
        Consequently the minimal profiles are geodesics for the conformal metric
    \begin{equation}\label{sf:metric}
        e^{-2V}(d\xi^2+dv^2),\qquad
            V=n\log\cosh\frac {\xi} {q}-a\log\sin v-b\log\cos v,
    \end{equation}
        up to an irrelevant constant factor.
        In Euclidean arclength $\ell$ and tangent angle $\alpha$,
            the profile equations \cite[\S2]{LW} become
    \begin{equation}\label{sf:ode}
        \dot\xi=\cos\alpha,\qquad \dot v=\sin\alpha,\qquad
            \dot\alpha=\mathfrak a(\xi)\sin\alpha+G(v)\cos\alpha.
    \end{equation}

    \begin{definition}\label{sf:admissible}
        An \emph{admissible half-arc} is a solution of \eqref{sf:ode} on a finite interval $[0,T]$ with $(\xi,v,\alpha)(0)=(0,\delta,0)$,
            $(\xi,v,\alpha)(T)=(0,\beta,\pi)$,
            where $0<\delta<\beta<\frac {\pi} {2}$ and $\xi,\dot v>0$ on $(0,T)$.
        Its union with its reflection across $\xi=0$ is the associated closed profile.
    \end{definition}

        By reflection symmetry and uniqueness,
            the closed profile is smooth.
        Since $\dot v>0$,
            we have $0<\alpha<\pi$ on the open half-arc,
            which is a graph $\xi=f(v)>0$ with $f'=\cot\alpha$,
            zero endpoint values and endpoint slopes $+\infty$ and $-\infty$.
        Thus the profile is simple,
            and the graph satisfies
    \begin{equation}\label{sf:graph}
        f''=-(1+f'^2)\bigl(\mathfrak a(f)+G(v)f'\bigr).
    \end{equation}

        For the shooting profile constructed in \cite[Proposition~3.4]{LW},
            the first return to $\xi=0$ coincides with the first positive zero of $\dot v$.
        Reflection therefore produces an admissible closed profile.

        Every admissible graph is strictly concave.
        Since $\mathfrak a'\leq M+1<Q$ by $Q\geq(\sqrt a+\sqrt b)^2\geq M+2$,
            differentiating \eqref{sf:ode} gives
    $$
        \frac {d} {d\ell}\bigl(e^{-V}\dot\alpha\bigr)
            =(\mathfrak a'-Q)e^{-V}\sin\alpha\cos\alpha.
    $$
        The quantity $e^{-V}\dot\alpha$ decreases near the initial endpoint and increases near the final one.
        At an interior minimum,
            $\cos\alpha=0$ and $\dot\alpha=\mathfrak a(\xi)>0$.
        Thus $\dot\alpha>0$ on the closed half-arc,
            so $f''<0$ and the reflected profile is strictly convex.
        At the endpoints,
            \eqref{sf:ode} gives $G(\delta)>0>G(\beta)$,
            so $\delta<v_*<\beta$.

        To control the lifted curvature,
            we bound the size of the profile and its distance from the chamber walls.
        Let $w=\max f$ and $h=\beta-\delta$.

    \begin{proposition}[Uniform profile estimates]\label{sf:profile-estimates}
        There are absolute constants $c,C>0$ such that every admissible profile satisfies
    $$
        \begin{gathered}
            c\leq\sqrt M\,w\leq C,\qquad
                \sqrt M\,h\leq C,\qquad
                \sup_{[\delta,\beta]}(G^2+Q)\leq CM,\\
            \delta\geq c\sqrt{\frac {a} {M}},\qquad
                \frac {\pi} {2}-\beta\geq c\sqrt{\frac {b} {M}}.
        \end{gathered}
    $$
    \end{proposition}
    \begin{proof}
        Let $\rho=\sin^a v\cos^b v$ and $\rho_*=\rho(v_*)$.
        By the arithmetic--geometric mean inequality,
    $$
        Q-G^2=M+2ab-a(a-1)\cot^2v-b(b-1)\tan^2v\leq3M.
    $$
        Thus $\rho''+3M\rho\geq0$.
        Integrating $(\rho'^2+3M\rho^2)'=2\rho'(\rho''+3M\rho)$ toward the maximum at $v_*$,
            we obtain $|\rho'|\leq\sqrt{3M}\rho_*$.
        Gauss--Bonnet and the identity
            $\frac {d} {d\ell}(\rho\cos\alpha)=-\rho\mathfrak a(\xi)\sin^2\alpha$
            give
    \begin{equation}\label{sf:gaussbonnet}
        \begin{gathered}
            \pi=\int_\delta^\beta(\mathfrak a(f)+Qf)\,dv,\\
            \rho(\delta)+\rho(\beta)
                =\int_\delta^\beta\frac {\rho\mathfrak a(f)} {\sqrt{1+f'^2}}\,dv.
        \end{gathered}
    \end{equation}
        Since $f$ is concave,
            $Q\geq M$ and $\mathfrak a(f)\leq(M+1)f$,
            we obtain
    $$
        Mwh\leq2\pi,
        \qquad
        \pi\geq\left(1+\frac {M} {M+1}\right)
            \frac {\rho(\delta)+\rho(\beta)} {\rho_*}.
    $$
        The derivative bound gives $\rho(\delta)+\rho(\beta)\geq(2-\sqrt{3M}h)\rho_*$,
            so $\sqrt{3M}h\geq2-\frac {3\pi} {5}>0$ for $M\geq2$.
        Consequently $w\sqrt M\leq C$ by $Mwh\leq2\pi$.
        Since $f$ is uniformly bounded and $q\in\{1,2\}$,
            we have $\mathfrak a(f)\geq c_0Mf$ for an absolute constant $c_0>0$.

        We next show that $\rho$ stays uniformly comparable to $\rho_*$ on the profile.
        Let $v_0=\frac {\delta+\beta} {2}$ and $T(v)=\min\{v-\delta,\beta-v\}$.
        Between $v_*$ and $v_0$,
            sine and cosine are at least half their values at $v_*$,
            so $Q\leq8M$ and $\rho(v_0)\geq e^{-Mh^2}\rho_*$.
        Using $f\geq\frac {w} {h}T$ and
            $\int_\delta^\beta QT\,dv=\log\frac {\rho(v_0)^2} {\rho(\delta)\rho(\beta)}$
            in \eqref{sf:gaussbonnet},
            we obtain
    $$
        Mw^2\log\frac {\rho_*^2} {\rho(\delta)\rho(\beta)}
            \leq2(Mwh)^2+\pi Mwh\leq10\pi^2.
    $$
        If either endpoint weight were $o(\rho_*)$ along a sequence,
            this would force $w\sqrt M\to0$.
        We have $\rho\leq\rho_*e^{-\frac {M} {2}(v-v_*)^2}$ since $Q\geq M$.
        Therefore,
            by the flux identity in \eqref{sf:gaussbonnet},
            $\rho(\delta)+\rho(\beta)\leq2\sqrt{2\pi M}\,w\rho_*=o(\rho_*)$.

        Concavity and the graph equation imply $-(\rho f')'\geq c_0M\rho f$.
        By Picone's identity and approximation,
            we have
    \begin{equation}\label{sf:gap}
        \int_\delta^\beta\rho\varphi'^2\,dv
            \geq c_0M\int_\delta^\beta\rho\varphi^2\,dv
            \qquad (\varphi\in H_0^1((\delta,\beta))).
    \end{equation}
        Along this sequence,
            let $L=\log\frac {\rho_*} {\max\{\rho(\delta),\rho(\beta)\}}\to\infty$
            and use the $H_0^1$ test function $\varphi=[1-\frac {1} {L}\log\frac {\rho_*} {\rho}]_+$ in \eqref{sf:gap}.
        By monotonicity on each side and the derivative bound,
            we have
    $$
        \int_\delta^\beta\rho\varphi'^2\,dv
            \leq\frac {2\sqrt{3M}\rho_*} {L}.
    $$
        For large $L$,
            the derivative bound gives $J=\{|v-v_*|\leq\frac {1} {2\sqrt{3M}}\}\subset(\delta,\beta)$,
            with $\rho\geq\frac {1} {2}\rho_*$ and $\varphi\geq\frac {1} {2}$ on $J$.
        Thus $\int_\delta^\beta\rho\varphi^2\,dv\geq c\frac {\rho_*} {\sqrt M}$,
            contradicting \eqref{sf:gap}.

        Consequently $\rho\geq c_1\rho_*$ throughout the interval.
        The derivative bound implies $G^2\leq CM$,
            while the two bounds on the distances to the chamber walls follow from $Q\leq G^2+3M\leq CM$.
        By the same flux bound,
            we have $w\sqrt M\geq c$.
        The upper bound on $h$ follows from $Mwh\leq2\pi$.
    \end{proof}

        For $U=\log W=-d\log\cosh(\frac {\xi} {q})+\log\rho(v)$,
            direct differentiation and Proposition~\ref{sf:profile-estimates} give $|\nabla U|_{g_Q}^2+|\Hess_{g_Q}U|_{g_Q}\leq CM$ along the profile.
        By stationarity,
            $|A_\Gamma|^2=(\eta U)^2$,
            and the curvature identity \eqref{mk:hessian-curvature},
            also valid for this suspension,
            yields $\sup_L|A_L|^2\leq CM$.
        Applying \eqref{eq:criterion} with $n=qM+1$ completes the proof of Theorem~\ref{thm:single} after increasing $M_0$.

\subsection{The extended families}

        For a fixed Gaussian seed in Theorem~\ref{thm:gaussian-seed},
            we have $\sup_L|A_L|^2\leq CD$,
            while $n+1=D+2k-1$.

        On the common compact annulus of \S\ref{sec:genus},
            the metrics $a_d^{-2}g_{\Sph^3}$ and $\widehat g_{\frac {2} {d}}$ have smooth positive limits.
        By Lemma~\ref{hg:fixed} and uniform comparison of these metrics and their connections,
            we have $\sup_{\Gamma_{d,h}}|A|^2_{g_{\Sph^3}}\leq C_h a_d^{-2}$.
        Since $u^2=1-a_d^2|x|^2$ and $da_d^2\leq2$,
            the spherical gradient satisfies $d|\nabla\log u|^2=d\frac {1-u^2} {u^2}=O_h(1)$.
        Thus \eqref{hg:shape} gives $\sup_{L_{d,h}}|A_L|^2\leq C_h d$,
            with $n+1=d+3$.

        The criterion \eqref{eq:criterion} now proves the area-minimizing assertions in Theorems~\ref{thm:families} and~\ref{thm:genus},
            after increasing $\mu_0(k)$ and $d_0(h)$,
            respectively.
    \begin{corollary}[Topological abundance in fixed dimension]\label{cor:abundance}
        For every fixed $k\geq2$ there are constants $c_k,D_k>0$ such that,
            for every integer $D\geq D_k$,
            $\R^{D+2k}$ contains at least $c_kD^{2k-1}$ area-minimizing hypercones,
            smooth away from the vertex,
            whose links are simply connected,
            inhomogeneous and of pairwise distinct homotopy types.
    \end{corollary}

    \begin{proof}
        For $g=2$,
            each input is a product of two spheres,
            so $k$ inputs contribute $r=2k$ sphere dimensions $m_1,\ldots,m_r$ with $\sum_jm_j=D$.
        Choose an integer $q\geq\max\{2,\mu_0(k)\}$.
        By the stars-and-bars formula,
            there are $\binom{D-rq+r-1}{r-1}$ ordered tuples with $m_j\geq q$ and $\sum_jm_j=D$.
        Each multiset has at most $r!$ orderings,
            so for sufficiently large $D$ the number of distinct multisets is at least
    $$
        \frac{1}{r!}\binom{D-rq+r-1}{r-1}\geq c_kD^{r-1}=c_kD^{2k-1}.
    $$
        Apply Theorem~\ref{thm:families}\textup{(a)} with Clifford leaves $P_i\cong\Sph^{m_{2i-1}}\times\Sph^{m_{2i}}$.
        The resulting links are $L\cong\Sph^{2k-2}\times\prod_{j=1}^{2k}\Sph^{m_j}$.
        The Poincar\'e polynomial of $L$ over $\mathbb Z_2$ is $(1+t^{2k-2})\prod_j(1+t^{m_j})\in\mathbb Z[t]$.
        After cancelling the common factor,
            we recover the multiset inductively from the smallest positive degree and its coefficient.
        Thus distinct multisets correspond to distinct homotopy types of links.
        The remaining assertions follow from $m_j\geq2$ and Remark~\ref{top:type}.
    \end{proof}

    \Needspace{11\baselineskip}
\subsection{Genus bounds at fixed dimension}

        In the construction of Theorem~\ref{thm:genus},
            we fix the genus before increasing the dimension.

    \begin{proposition}\label{prop:fixed-genus}
        Fix $d\geq1$,
            and let $\Gamma\subset Q=\{u>0\}\subset\Sph^3$ be a smooth closed connected embedded surface of genus $h$.
        If the original oriented cone over its lift \eqref{hg:lift} is area minimizing,
            then,
            with $N=d+3$ and $\omega_j=\operatorname{Area}(\Sph^j)$,
    $$
        \int_\Gamma u^d\,dA\leq\frac {N\omega_N} {2\omega_d},\qquad
            (h-1)\bigl(\min_\Gamma u\bigr)^d
            \leq\frac {N(N-2)^2\omega_N} {64\pi\omega_d}.
    $$
        Consequently the genus is uniformly bounded for such surfaces contained in any fixed compact set $K\Subset Q$,
            whereas any such fixed-$d$ sequence with $h\to\infty$ must satisfy $\min_\Gamma u\to0$.
    \end{proposition}

    \begin{proof}
        Since the lift $L\subset\Sph^N$ separates the sphere,
            the smaller complementary region,
            competes with $C(L)\cap B_1$ and has mass at most $\frac {\omega_N} {2}$.
        Hence $\frac {\omega_d} {N}\int_\Gamma u^d\,dA=\frac {\operatorname{Area}(L)} {N}\leq\frac {\omega_N} {2}$
            which leads to the first bound.
        Let $u_0=\min_\Gamma u>0$.
        Then, according to (1) Gauss equation and Gauss--Bonnet,
            together with \eqref{hg:shape} and the fiber volume formula,
           and (2) the cone stability and the radial Hardy inequality applied to the constant link function 
          we get 
    $$
        8\pi(h-1)u_0^d\omega_d
            \leq\int_L|A_L|^2\,dA_L
            \leq\frac {(N-2)^2} {4}\operatorname{Area}(L).
    $$
        Combining this with the volume bound proves the genus estimate.
    \end{proof}

    \Needspace{8\baselineskip}
    
\end{document}